\documentclass[10pt,leqno]{amsart}
\usepackage{amsfonts}
\usepackage{amsmath}
\usepackage{amsthm}
\usepackage{amssymb}

\theoremstyle{definition}
\newtheorem{definition}{Definition}[section]

\theoremstyle{plain}
\newtheorem{thm}[definition]{Theorem}

\newtheorem{pro}[definition]{Proposition}
\newtheorem{cor}[definition]{Corollary}

\begin{document}
\def\iff{if and only if }
\def\pi{positive implicative }

\setcounter{page}{1} 
\vspace{16mm}

\begin{center}
{\normalsize \textbf{Euclidean Distance between Two Linear Varieties}} \\[%
12mm]
\textsc{M. A. Facas Vicente$^{\ast }$$^{1}$$%
^{2}$, Armando Gon\c calves$^{1}$ and Jos\'{e} Vit\'{o}ria$^{1}$}\\[8mm]

\begin{minipage}{123mm}
{\small {\sc Abstract.}
   This paper addresses to the problem of finding the (minimum) Euclidean distance between two linear varieties. This problem is, usually, solved minimizing a target function. We propose a novel approach: to use the Moore-Penrose generalized inverse of a matrix in order to find the Euclidean distance and, at the same time, going further by exhibiting the best approximation pair of points.}
\end{minipage}
\end{center}

\bigskip

\bigskip

\renewcommand{\thefootnote}{} \footnotetext{$^{\ast }$\thinspace
Corresponding author.} \footnotetext{%
2000 \textit{Mathematics Subject Classification.} 41A17, 51M16.}
\footnotetext{\textit{Key words and phrases.} linear varieties, Euclidean distance, Moore-Penrose inverse, best
approximation pair.} \footnotetext{$^{1}$\thinspace
Department of Mathematics, University of Coimbra, Apartado 3008, EC\ Santa
Cruz, 3001-501 Coimbra, Portugal. E-mails: \textit{vicente@mat.uc.pt} (M. A. Facas Vicente), \textit{adsg@mat.uc.pt} (Armando Gon\c calves) and \textit{jvitoria@mat.uc.pt} (Jos\'{e} Vit\'{o}ria).}
\footnotetext{$^{2}$\thinspace Supported by Instituto de Engenharia de
Sistemas e Computadores---Coimbra, Rua Antero de Quental, 199, 3000-033
Coimbra, Portugal.}


\section{Introduction}


\bigskip

The distance between two linear varieties has been dealt with in several papers: \cite{Caseiro}, \cite{Dax} and \cite{LAA1}. In \cite{LAA1}, only the distance is considered and the points that realize the distance are not exhibited. In \cite{Caseiro} the best approximation points are found.

In this paper we present a novel method to obtain the Euclidean distance between two linear varieties. We go further than the usual: we exhibit the two points, on both varieties, that materialize the minimum distance. That our method may have several computational advantages is worth to be treated elsewhere.

The problem of finding the best approximation pair of points between two linear varieties $V_{\overrightarrow{b}}$ and $V_{\overrightarrow{c}}$, is solved by using the Moore-Penrose generalized inverse of a matrix.

\begin{definition} The Moore-Penrose inverse $M^\dag$ of a real matrix $M$ is the unique matrix that satisfies the following four equations:
\begin{enumerate}
\item $MM^\dag M=M$;
\item $M^\dag M M^\dag =M^\dag $;
\item $(M ^\dag )^T =M M^\dag $;
\item $(M^\dag M)^T = M^\dag M$,
\end{enumerate}
where $T$ stands for the matrix transposition.
\end{definition}

From \cite[pages 28, 29 and 58]{Campbell} we can extract the following result.

\begin{thm}
Let be given the matrix $A\in \mathrm{I\kern-.17emR}^{m\times n}$ and the vector $\overrightarrow{b} \in \mathrm{I\kern-.17emR}^m $. We have the following:
\begin{enumerate}
\item $A^\dag \overrightarrow{b} $ is the minimal least squares solution to the system $A\overrightarrow{x}=\overrightarrow{b} $;
\item The vector $\overrightarrow{u}$ is a least squares solution of $A\overrightarrow{x}=\overrightarrow{b} $ if and only if $\overrightarrow{u} $ is of the form $A^\dag \overrightarrow{b} + \overrightarrow{h} $, where $\overrightarrow{h}\in \mbox{ker}(A)$;
\item If $C\in \mathrm{I\kern-.17emR}^{m\times p}$, the Moore-Penrose inverse of $\left[A\,\vdots \, C\right]$ can be written in the following way
    $$\left[A\,\vdots \, C\right]^\dag =
\!\left[\begin{array}{c}
(I+L L^T)^{-1} (A^\dag - A^\dag C B^\dag )\\
\hline \\
L^T (I+L L^T)^{-1} (A^\dag - A^\dag C B^\dag )+ B^\dag
\end{array}\right]\\,
$$
where $B=(I-AA^\dag )C$ and $L=A^\dag C(I-B^\dag B).$
\end{enumerate}
\end{thm}

Then the Moore-Penrose generalized inverse gives the minimal least squares solution of a linear system of equations (see also \cite{Greville} and \cite{Meyer}).

Furthermore, as the matrix of the system to be solved is a column partitioned matrix formed by the matrices defining the linear varieties in the way of \cite{Dax}, we can use Moore-Penrose inverses of smaller matrices to find the Moore-Penrose inverse of the system matrix. In this manner, some computational advantages are expected to be obtained due to the fact that computing the Moore-Penrose inverse is difficult because, on one side it does not need to be a continuous function of the entries of the matrix \cite[pages 423-424]{Meyer}, on the other it is not numerically stable \cite[page 247]{Campbell} and \cite[pages 423-424]{Meyer}.

\bigskip

\bigskip

\section{Euclidean Distance between two Linear Varieties}

\bigskip

Let the varieties $V_{\overrightarrow{b}}$ and $V_{\overrightarrow{c}}$ in $\mathrm{I\kern-.17emR}^{m}$ be given by
\begin{equation*}
V_{\overrightarrow{b}}:=\left\{ \overrightarrow{b}+B\overrightarrow{u}:%
\overrightarrow{u}\in \mathrm{I\kern-.17emR}^{l_{1}}\right\}
\end{equation*}
and
\begin{equation*}
V_{\overrightarrow{c}}:=\left\{ \overrightarrow{c}-C\overrightarrow{v}:%
\overrightarrow{v}\in \mathrm{I\kern-.17emR}^{l_{2}}\right\},
\end{equation*}
where $\overrightarrow{b}$ and $\overrightarrow{c}$ are given vectors in $\mathrm{I\kern-.17emR}^{m}$ and $B$ and $C$ are real matrices of type $m\times l_{1}$ and $m\times l_{2}$, respectively.

The Euclidean distance $d\left( V_{\overrightarrow{b}},V_{\overrightarrow{c}}\right) $ is obtained by solving the following problem  \cite{Dax}:
\begin{equation}
\text{\texttt{minimize } }f\left( \overrightarrow{u},\overrightarrow{v}\right) =\left\Vert
\overrightarrow{b}+B\overrightarrow{u}-\overrightarrow{c}+C\overrightarrow{v}%
\right\Vert .  \label{eq1}
\end{equation}

\begin{definition}
The best approximation pair of the linear varieties $V_{\overrightarrow{b}}$ and $V_{\overrightarrow{c}}$ is the pair $(\overrightarrow{b}^{\ast } , \overrightarrow{c}^{\ast }) \in V_{\overrightarrow{b}} \times V_{\overrightarrow{c}}$ such that $$\| \overrightarrow{b}^{\ast } - \overrightarrow{c}^{\ast } \| = d(V_{\overrightarrow{b}} , V_{\overrightarrow{c}}).$$
\end{definition}

Problem (1) is equivalent to the problem
\begin{equation}
\text{\texttt{minimize} }\left\Vert A\overrightarrow{x}-\overrightarrow{d}%
\right\Vert ,  \label{eq2}
\end{equation}
where
\begin{equation*}
\overrightarrow{d}=\overrightarrow{c}-\overrightarrow{b}
\end{equation*}
and
\begin{equation*}
A=\left[  B\, \vdots \, C \right]
\end{equation*}
is a $m\times n$ real matrix, $n=l_{1}+l_{2}$ and
\begin{equation*}
\overrightarrow{x}=\left[
\begin{array}{c}
\overrightarrow{u} \\
\overrightarrow{v}%
\end{array}%
\right] \in \mathrm{I\kern-.17emR}^{n}.
\end{equation*}

Problem (\ref{eq2}) is equivalent to find the least norm solution $\overrightarrow{x}^{\ast }$ of the linear system of equations
\begin{equation*}
A\overrightarrow{x}=\overrightarrow{d}.  \label{eq3}
\end{equation*}

By using the Moore-Penrose inverse, we get \cite[page 109]{Greville}, \cite[page 28]{Campbell}, \cite[page 438]{Meyer}
\begin{equation*}
\overrightarrow{x}^{\ast }=A^{\dag }\overrightarrow{d},  \label{eq4}
\end{equation*}%
with
\begin{equation*}
\overrightarrow{x}^{\ast }=\left[
\begin{array}{c}
\overrightarrow{u}^{\ast } \\
\overrightarrow{v}^{\ast }%
\end{array}%
\right] .
\end{equation*}%

The above developments allow the formalization of our main result.

\begin{pro}
Let be given two skew linear varieties
\begin{equation*}
V_{\overrightarrow{b}}:=\left\{ \overrightarrow{b}+B\overrightarrow{u}:%
\overrightarrow{u}\in \mathrm{I\kern-.17emR}^{l_{1}}\right\}
\end{equation*}
and
\begin{equation*}
V_{\overrightarrow{c}}:=\left\{ \overrightarrow{c}-C\overrightarrow{v}:%
\overrightarrow{v}\in \mathrm{I\kern-.17emR}^{l_{2}}\right\},
\end{equation*}
where $\overrightarrow{b}$ and $\overrightarrow{c}$ are given vectors in $\mathrm{I\kern-.17emR}^{m}$ and $B$ and $C$ are real matrices of type $m\times l_{1}$ and $m\times l_{2}$, respectively.

Then, the best approximation pair  $(\overrightarrow{b}^{\ast } , \overrightarrow{c}^{\ast }) \in V_{\overrightarrow{b}} \times V_{\overrightarrow{c}}$ is given by
$$ \overrightarrow{b}^{\ast }=\overrightarrow{b}+B\overrightarrow{u}^{\ast }\ and \,\, \overrightarrow{c}^{\ast }=\overrightarrow{c}-C\overrightarrow{v}^{\ast },$$
where \begin{equation*}
\left[
\begin{array}{c}
\overrightarrow{u}^{\ast } \\
\overrightarrow{v}^{\ast }
\end{array}
\right]:= \overrightarrow{x}^{\ast }=A^\dag \overrightarrow{d},
\end{equation*}
with $A=\left[B\,\vdots \, C \right] \in \mathrm{I\kern-.17emR}^{m\times (l_1 + l_2 )}$ and $\overrightarrow{d}=\overrightarrow{c} - \overrightarrow{b} $.

Furthermore, the distance between $V_{\overrightarrow{b}}$ and $V_{\overrightarrow{c}}$ is $d(V_{\overrightarrow{b}},\,V_{\overrightarrow{c}} ) = \| \overrightarrow{b}^{\ast } - \overrightarrow{c}^{\ast } \|$.
\end{pro}

In order to make the most of the partitioned structure of the matrix $A=\left[B\,\vdots \, C\right] $ and in hopes of bettering the computational aspects, we have the following result.

\begin{cor}
Let be given the non-intersecting and non-parallel linear varieties

\begin{equation*}
V_{\overrightarrow{b}}:=\left\{ \overrightarrow{b}+B\overrightarrow{u}:%
\overrightarrow{u}\in \mathrm{I\kern-.17emR}^{l_{1}}\right\}
\end{equation*}
and
\begin{equation*}
V_{\overrightarrow{c}}:=\left\{ \overrightarrow{c}-C\overrightarrow{v}:%
\overrightarrow{v}\in \mathrm{I\kern-.17emR}^{l_{2}}\right\},
\end{equation*}
where $\overrightarrow{b}$ and $\overrightarrow{c}$ are given vectors in $\mathrm{I\kern-.17emR}^{m}$ and $B$ and $C$ are real matrices of type $m\times l_{1}$ and $m\times l_{2}$, respectively.

Then, the best approximation points $\overrightarrow{b}^{\ast } = \overrightarrow{b} + B \overrightarrow{u}^{\ast}$ and $\overrightarrow{c}^{\ast } = \overrightarrow{c} - C \overrightarrow{v}^{\ast} $ are obtained through
$$ \left[
\begin{array}{c}
\overrightarrow{u}^{\ast } \\
\overrightarrow{v}^{\ast }
\end{array}
\right] = \left[\begin{array}{c}
(I+H H^T)^{-1} (B^\dag - B^\dag C G^\dag )\\
\hline \\
H^T (I+H H^T)^{-1} (B^\dag - B^\dag C G^\dag )+ G^\dag
\end{array}\right] (\overrightarrow{c} - \overrightarrow{b}),$$
with $G=(I-BB^\dag )C$ and $H=B^\dag C(I-G^\dag G).$
\end{cor}

Notice that the distance between the two linear varieties is not given by
$ \|\overrightarrow{x}^{\ast }\|$, as we can see if we consider, for example, the straight lines defined, in $\mathrm{I\kern-.17emR}^3 $, by $y=0, z=1$ and $x=0, z=0$.

\bigskip

\section{Conclusions}

Use was made of the Moore-Penrose inverse of a matrix to find the distance between two skew linear varieties in $\mathrm{I\kern-.17emR}^m $. This novel approach was made possible by a formulation of the distance problem by considering the difference set of two closed convex sets.

Computational difficulties were mentioned. However, some moderate computational gains can be expected, due to the fact that use is made of a formula for the Moore-Penrose inverse of a partitioned matrix.

\end{document}